\documentclass[12pt,oneside]{article}
\usepackage{amsmath,amssymb,amsfonts,amsthm}
\newcommand{\T}{{\cal T}}

\newcommand{\Real}{\mathbb R}

\newcommand{\set}[1]{\left\{#1\right\}}

\newcommand{\tm}{\T M}
\newcommand{\p}{\pi^{-1}(TM)}
\newcommand {\cp}{\mathfrak{X}(\pi (M))}

\def\o#1{\overline{#1}}

\def\Section#1{\vspace{30truept}\addtocounter{section}{1}\setcounter{thm}{0}
\setcounter{equation}{0}{\noindent\Large\bf
    \arabic{section}.~~#1}\par \vspace{12pt}}

\newtheorem{thm}{Theorem}[section]
\newtheorem{cor}[thm]{Corollary}
\newtheorem{lem}[thm]{Lemma}

\numberwithin{equation}{section}

\begin{document}
\title{{\textbf{On Akbar-Zadeh's Theorem on  a Finsler  Space of Constant Curvature}}} 
\author{\bf A. Soleiman}
\date{}

\maketitle                     
\vspace{-1cm}
\begin{center}
{Department of Mathematics, Faculty of Science,\\ Benha
University, Benha,
 Egypt}
\end{center}
\vspace{-0.8cm}

\begin{center}
soleiman@yahoo.com
\end{center}
\smallskip
\vspace{1cm} \maketitle

\smallskip
\bigskip

\noindent{\bf Abstract.}  The aim of the present paper is to give
two \emph{intrinsic} generalizations of Akbar-Zadeh's theorem on a
Finsler space of constant curvature. Some consequences, of these
generalizations, are drown.



\vspace{30truept}\centerline{\Large\bf{Introduction}}\vspace{12pt}
\par

 In \cite{r58a}, Akbar-Zadeh proved \emph{locally} that if the  $h$-curvature $R^{r}_{ijk}$ of the Cartan
connection $C\Gamma$ associated with a Finsler manifold $(M,L)$,
$dim M\geq3$, satisfies
$$R^{r}_{ijk}=k(x,y)(g_{ij}\delta^{r}_{k}-g_{ik}\delta^{r}_{j}),$$
 where $k(x,y)$ is a  scalar function on $\tm$, positively homogeneous of degree zero ((0) p-homogeneous), then
\begin{description}
\item[(a)] $k$ is constant,

\item[(b)] if $k\neq0$, then
 \begin{description}
  \item[(1)] the $v$-curvature of $C\Gamma$ vanishes: $S^{r}_{ijk}=0,$
  \item[(2)] the $hv$-curvature of $C\Gamma$ is symmetric with respect to the  last two indices:  $P^{r}_{ijk}=P^{r}_{ikj}$.
\end{description}
\end{description}
\par
In \cite{hojo}, H\={o}j\={o} showed \emph{locally} that if the
$h$-curvature $R^{r}_{ijk}$ of the  generalized Cartan connection
$C\Gamma$, $dim M\geq3$, satisfies\footnote{$\mathfrak{A}_{ij}$
indicates interchanges of indices $j$ and $k$, and subtraction:
$\mathfrak{A}_{ij} \set{F_{ij}}=F_{ij}-F_{ji}$}
$$R^{r}_{ijk}=k(x,y)\mathfrak{A}_{j,k}\set{q\,g_{ij}\delta^{r}_{k}+(q-2)(g_{ij}\ell_{k}\ell^{r}-\delta^{r}_{j}\ell_{i}\ell_{k}},$$
 where $k(x,y)$ is a (0) p-homogeneous scalar function and $1\neq q\in\Real$, then
\begin{description}
\item[(a)] $k$ is constant,

\item[(b)] if $k\neq0$, then
 \begin{description}
  \item[(1)] the $v$-curvature of $C\Gamma$ satisfies  $S^{r}_{ijk}=\frac{q-2}{2(1-q)}\mathfrak{A}_{j,k}\set{\hbar_{ij}\hbar^{r}_{k} },$
  \item[(2)] the $hv$-curvature of $C\Gamma$ is symmetric with respect to the  last two indices.
\end{description}
\end{description}
\par
The aim of the present paper is to provide two \emph{intrinsic}
generalizations of Akbar-Zadeh's and H\={o}j\={o}'s theorems. As a
by-product, some  consequences concerning $S_{3}$-like and
$S_{4}$-like spaces, are deduced.
\par
The present work is formulated in a coordinate-free form, without
being trapped into the complications of indices. However, the local
expressions of the obtained results, when calculated, coincides with
the existing local results.

\Section{Notation and Preliminaries}

In this section, we give a brief account of the basic concepts
 of the pullback approach to intrinsic Finsler geometry necessary for this work. For more
 details, we refer to \cite{r58}, \cite{r61} and \cite{amr}. We
 shall use the same notations of \textbf{\cite{amr}}.

 In what follows, we denote by $\pi: \T M\longrightarrow M$ the subbundle of nonzero vectors
tangent to $M$, $\cp$
  the $\mathfrak{F}(TM)$-module of differentiable sections of the pullback bundle $\pi^{-1}(T M)$.
The elements of $\mathfrak{X}(\pi (M))$ will be called $\pi$-vector
fields and will be denoted by barred letters $\overline{X} $. The
tensor fields on $\pi^{-1}(TM)$ will be called $\pi$-tensor fields.
The fundamental $\pi$-vector field is the $\pi$-vector field
$\overline{\eta}$ defined by $\overline{\eta}(u)=(u,u)$ for all
$u\in TM$.
\par
We have the following short exact sequence of vector bundles
$$0\longrightarrow
 \pi^{-1}(TM)\stackrel{\gamma}\longrightarrow T(\T M)\stackrel{\rho}\longrightarrow
\pi^{-1}(TM)\longrightarrow 0 ,\vspace{-0.1cm}$$ with the well known
definitions of  the bundle morphisms $\rho$ and $\gamma$. The vector
space $V_u (\T M)= \{ X \in T_u (\T M) : d\pi(X)=0 \}$  is called
the vertical space to $M$ at $u$.
\par
Let $D$ be  a linear connection (or simply a connection) on the
pullback bundle $\pi^{-1}(TM)$.
 We associate with
$D$ the map \vspace{-0.1cm} $K:T \T M\longrightarrow
\pi^{-1}(TM):X\longmapsto D_X \overline{\eta} ,$ called the
connection map of $D$.  The vector space $H_u (\T M)= \{ X \in T_u
(\T M) : K(X)=0 \}$ is called the horizontal space to $M$ at $u$ .
   The connection $D$ is said to be regular if
$ T_u (\T M)=V_u (\T M)\oplus H_u (\T M) \,  \forall \, u\in \T M$.
\par If $M$ is endowed with a regular connection, then the vector bundle
   maps $
 \gamma,
   \rho |_{H(\T M)}$ and $
   K |_{V(\T M)}$
 are vector bundle isomorphisms. The map
 $\beta:=(\rho |_{H(\T M)})^{-1}$
 will be called the horizontal map of the connection
$D$.
\par
 The horizontal ((h)h-) and
mixed ((h)hv-) torsion tensors of $D$, denoted by $Q $ and $ T $
respectively, are defined by \vspace{-0.2cm}
$$Q (\overline{X},\overline{Y})=\textbf{T}(\beta \overline{X}\beta \overline{Y}),
\, \,\, T(\overline{X},\overline{Y})=\textbf{T}(\gamma
\overline{X},\beta \overline{Y}) \quad \forall \,
\overline{X},\overline{Y}\in\mathfrak{X} (\pi (M)),\vspace{-0.2cm}$$
where $\textbf{T}$ is the (classical) torsion tensor field
associated with $D$.
\par
The horizontal (h-), mixed (hv-) and vertical (v-) curvature tensors
of $D$, denoted by $R$, $P$ and $S$ respectively, are defined by
$$R(\overline{X},\overline{Y})\overline{Z}=\textbf{K}(\beta
\overline{X}\beta \overline{Y})\overline{Z},\quad
P(\overline{X},\overline{Y})\overline{Z}=\textbf{K}(\beta
\overline{X},\gamma \overline{Y})\overline{Z},\quad
S(\overline{X},\overline{Y})\overline{Z}=\textbf{K}(\gamma
\overline{X},\gamma \overline{Y})\overline{Z}, $$
 where $\textbf{K}$
is the (classical) curvature tensor field associated with $D$.
\par
The contracted curvature tensors of $D$, denoted by $\widehat{R}$,
$\widehat{P}$ and $\widehat{S}$ respectively, known also as the
 (v)h-, (v)hv- and (v)v-torsion tensors, are defined by
$$\widehat{R}(\overline{X},\overline{Y})={R}(\overline{X},\overline{Y})\overline{\eta},\quad
\widehat{P}(\overline{X},\overline{Y})={P}(\overline{X},\overline{Y})\overline{\eta},\quad
\widehat{S}(\overline{X},\overline{Y})={S}(\overline{X},\overline{Y})\overline{\eta}.$$
\par
If $M$ is endowed with a metric $g$ on $\p$, we write
\begin{equation}\label{cur.g}
    R(\overline{X},\overline{Y},\overline{Z}, \overline{W}):
=g(R(\overline{X},\overline{Y})\overline{Z}, \overline{W}),\,
\cdots, \, S(\overline{X},\overline{Y},\overline{Z}, \overline{W}):
=g(S(\overline{X},\overline{Y})\overline{Z}, \overline{W}).
\end{equation}
\par
The following result is of extreme importance. \vspace{-0.1cm}
\begin{thm} {\em\cite{r94}} \label{th.1} Let $(M,L)$ be a Finsler
manifold and  $g$ the Finsler metric defined by $L$. There exists a
unique regular connection $\nabla$ on $\pi^{-1}(TM)$ such
that\vspace{-0.2cm}
\begin{description}
  \item[(a)]  $\nabla$ is  metric\,{\em:} $\nabla g=0$,

  \item[(b)] The (h)h-torsion of $\nabla$ vanishes\,{\em:} $Q=0
  $,
  \item[(c)] The (h)hv-torsion $T$ of $\nabla$\, satisfies\,\emph{:}
   $g(T(\overline{X},\overline{Y}), \overline{Z})=g(T(\overline{X},\overline{Z}),\overline{Y})$.
\end{description}
\par
 Such a connection is called the Cartan
connection associated with  the Finsler manifold $(M,L)$.
\end{thm}


\Section{First generalization of Akbar-Zadeh theorem}

 In this section, we investigate an intrinsic  generalization of Akbar-Zadeh
 theorem. We begin first with the following  two lemmas which will be useful
for subsequent use.

\begin{lem}\label{ricc} Let $\nabla$ be the Cartan connection on a Finsler manifold $(M,L)$. For a $\pi$-tensor field
$\omega$ of type $(1,1)$, we have the following commutation
formulae{\em:}\,
\begin{description}
 \item[\em{\textbf{(a)}}]
$(\stackrel{2}{\nabla}\stackrel{2}{\nabla}\omega)(\overline{X},\overline{Y},\overline{Z})-
(\stackrel{2}{\nabla}\stackrel{2}{\nabla}\omega)(\overline{Y},\overline{X},\overline{Z})=
\omega(S(\overline{X},\overline{Y})\overline{Z})-S(\overline{X},\overline{Y})\omega(\overline{Z}),$

\item[\em{\textbf{(b)}}]$(\stackrel{2}{\nabla}\stackrel{1}{\nabla}\omega)(\overline{X},\overline{Y},\overline{Z})-
(\stackrel{1}{\nabla}\stackrel{2}{\nabla}\omega)(\overline{Y},\overline{X},\overline{Z})=
\omega(P(\overline{X},\overline{Y})\overline{Z})-P(\overline{X},\overline{Y})\omega(\overline{Z})$\\
${\qquad\qquad\qquad\qquad\qquad\qquad\qquad\!\qquad}+
(\stackrel{2}{\nabla}\omega)(\widehat{P}(\overline{X},\overline{Y}),\overline{Z})+
(\stackrel{1}{\nabla}\omega)(T(\overline{Y},\overline{X}),\overline{Z})$,

 \item[\em{\textbf{(c)}}]$(\stackrel{1}{\nabla}\stackrel{1}{\nabla}\omega)(\overline{X},\overline{Y},\overline{Z})-
(\stackrel{1}{\nabla}\stackrel{1}{\nabla}\omega)(\overline{Y},\overline{X},\overline{Z})=
\omega(R(\overline{X},\overline{Y})\overline{Z})-R(\overline{X},\overline{Y})\omega(\overline{Z})$\\
${\qquad\qquad\qquad\qquad\qquad\qquad\qquad\!\qquad}+
(\stackrel{2}{\nabla}\omega)(\widehat{R}(\overline{X},\overline{Y}),\overline{Z})$,
\end{description}
where $\stackrel{1}\nabla$ and $\stackrel{2}\nabla$ are  the h- and
v-covariant derivatives associated with $\nabla$.
\end{lem}

\begin{lem}\label{ell}
Let $(M,L)$ be a Finsler manifold, $g$ the Finsler metric defined by
$L$, $\ell:=L^{-1}i_{\overline{\eta}}g$ and $\hbar:=\ell \circ
\ell-g$ the angular metric tensor. Then, we have:
\begin{description}
 \item[\em{\textbf{(a)}}]
$\stackrel{1}{\nabla}L=0,\quad\stackrel{2}{\nabla}L=\ell.$

 \item[\em{\textbf{(b)}}]
$\stackrel{1}{\nabla}\ell=0,
\quad\stackrel{2}{\nabla}\ell=L^{-1}\hbar.$

 \item[\em{\textbf{(c)}}]
$i_{\o \eta}\ell=L, \quad i_{\o \eta}\hbar=0.$

\end{description}
\end{lem}

\begin{proof} The proof follows from the fact that $\nabla g=0$ and $g(\o \eta,\o \eta)=L^{2}$.
\end{proof}

Now, we have
 \begin{thm}\label{prop.1} Let $(M,L)$ be a Finsler manifold of dimension $n$ and  $g$ the Finsler metric defined by $L$.
 If the $(v)h$-torsion tensor  $\widehat{R}$ of the Cartan connection  is given
 by
 \begin{equation}\label{eq.1}
    \widehat{R}(\overline{X},\overline{Y})=kL(\ell(\overline{X})\overline{Y}-\ell(\overline{Y})\overline{X}),
 \end{equation}
where $k(x,y)$ is  a homogenous function of degree $0$ in $y$
\emph{(}$\nabla_{\gamma \o \eta}k=0$\emph{)}, then we have:
\begin{description}
  \item[(a)] $\mathfrak{S}_{\overline{X},\overline{Y},\overline{Z}}R(\overline{X},\overline{Y})\overline{Z}=0.$
  \footnote{$\mathfrak{S}_{\overline{X},\overline{Y},\overline{Z}}$ denotes the cyclic sum over
   $\overline{X},\overline{Y}$ and $\overline{Z}$.}
  \item[(b)] $k$ is  constant if $\emph{dim}\,M\geq3$.
\end{description}
 \end{thm}

\begin{proof}~\par
\noindent\textbf \  {(a)} We have \cite{r96}:
\begin{equation}\label{eq.2}
   \mathfrak{S}_{\overline{X},\overline{Y},\overline{Z}}R(\overline{X},\overline{Y})\overline{Z}=
    \mathfrak{S}_{\overline{X},\overline{Y},\overline{Z}}T(\widehat{R}(\overline{X},\overline{Y}),\overline{Z}).
 \end{equation}
From (\ref{eq.1}), noting that the $(h)hv$-torsion $T$ is symmetric,
we obtain
\begin{eqnarray}
   \mathfrak{S}_{\overline{X},\overline{Y},\overline{Z}}T(\widehat{R}(\overline{X},\overline{Y}),\overline{Z})
    &=&kL T(\ell(\overline{X})\overline{Y}-\ell(\overline{Y})\overline{X},\overline{Z}) \nonumber\\
   && +kL T(\ell(\overline{Y})\overline{Z}-\ell(\overline{Z})\overline{Y},\overline{X})\nonumber \\
   && +kL
   T(\ell(\overline{Z})\overline{X}-\ell(\overline{X})\overline{Z},\overline{Y})\nonumber\\
   &=&kL\set{
   \ell(\overline{X})T(\overline{Y},\overline{Z})-
    \ell(\overline{Y})T(\overline{X},\overline{Z})}\nonumber\\
    &&+kL\set{
   \ell(\overline{Y})T(\overline{Z},\overline{X})-
    \ell(\overline{Z})T(\overline{Y},\overline{X})}\nonumber\\
    &&+kL\set{
   \ell(\overline{Z})T(\overline{X},\overline{Y})-
    \ell(\overline{X})T(\overline{Z},\overline{Y})}\nonumber\\
    &=&0. \label{eq.3}
\end{eqnarray}
Hence, the result follows from (\ref{eq.2}) and (\ref{eq.3}).

\vspace{5pt}
 \noindent\textbf{(b)} We have \cite{r96}:

\begin{equation}\label{eq.4}
   \mathfrak{S}_{\overline{X},\overline{Y},\overline{Z}}\,
\{(\nabla_{\beta \overline{X}}R)(\overline{Y},
\overline{Z},\overline{W})+P(\widehat{R}(\overline{X},\overline{Y}),
\overline{Z})\overline{W}\}=0.
\end{equation}
 From (\ref{eq.1}), noting that the $(v)hv$-torsion $\widehat{P}$ is symmetric \cite{r96}, we get
\begin{eqnarray*}
  \mathfrak{S}_{\overline{X},\overline{Y},\overline{Z}}\widehat{P}(\widehat{R}(\overline{X},\overline{Y}),
\overline{Z}) &=&kL\set{
   \ell(\overline{X})\widehat{P}(\overline{Y},\overline{Z})-
    \ell(\overline{Y})\widehat{P}(\overline{X},\overline{Z})}\nonumber\\
    &&+kL\set{
   \ell(\overline{Y})\widehat{P}(\overline{Z},\overline{X})-
    \ell(\overline{Z})\widehat{P}(\overline{Y},\overline{X})}\nonumber\\
    &&+kL\set{
   \ell(\overline{Z})\widehat{P}(\overline{X},\overline{Y})-
    \ell(\overline{X})\widehat{P}(\overline{Z},\overline{Y})}\nonumber\\
    &=&0. \label{eq.5}
\end{eqnarray*}
From which, together with (\ref{eq.4}), it follows that
\begin{equation}\label{eq.6}
   \mathfrak{S}_{\overline{X},\overline{Y},\overline{Z}}\,
(\nabla_{\beta \overline{X}}\widehat{R})(\overline{Y},
\overline{Z})=0.
\end{equation}
Again from (\ref{eq.1}), noting that $\nabla_{\beta
\overline{X}}\ell=0$ (Lemma \ref{ell}(b)), (\ref{eq.6}) reads
\begin{eqnarray*}
  && L(\nabla_{\beta
\overline{X}}k)\set{\ell(\overline{Y})\overline{Z}-\ell(\overline{Z})\overline{Y}}
+  L(\nabla_{\beta
\overline{Y}}k)\set{\ell(\overline{Z})\overline{X}-\ell(\overline{X})\overline{Z}} \\
   &&+L(\nabla_{\beta
\overline{Z}}k)\set{\ell(\overline{X})\overline{Y}-\ell(\overline{Y})\overline{X}}=0.
\end{eqnarray*}
Setting $\overline{Z}=\overline{\eta}$ into the above equation,
noting that $\ell(\overline{\eta})=L$  ((Lemma \ref{ell})(c)), we
obtain
\begin{eqnarray*}
  && L(\nabla_{\beta
\overline{X}}k)\set{\ell(\overline{Y})\overline{\eta}-L\overline{Y}}
+  L(\nabla_{\beta
\overline{Y}}k)\set{L\overline{X}-\ell(\overline{X})\overline{\eta}} \\
   &&+L(\nabla_{\beta
\overline{\eta}}k)\set{\ell(\overline{X})\overline{Y}-\ell(\overline{Y})\overline{X}}=0.
\end{eqnarray*}
Taking the trace of both sides with respect to $\overline{Y}$, it
follows that
\begin{equation}\label{eq.7}
\nabla_{\beta \overline{X}}k=L^{-1}(\nabla_{\beta
\overline{\eta}}k)\ell(\overline{X}).
\end{equation}
\par
On the other hand, we have \cite{r96}
\begin{eqnarray}\label{eq.r1a}
 &&(\nabla_{\gamma\overline{X}}R)(\overline{Y},\overline{Z},\overline{W})
   + (\nabla_{\beta\overline{Y}}P)(\overline{Z},\overline{X},\overline{W})-
   (\nabla_{\beta
   \overline{Z}}P)(\overline{Y},\overline{X},\overline{W})\nonumber\\
&&-
P(\overline{Z},\widehat{P}(\overline{Y},\overline{X}))\overline{W}
+R(T(\overline{X},\overline{Y}),\overline{Z})\overline{W}-
S(\widehat{R}(\overline{Y},\overline{Z}),\overline{X})\overline{W}\nonumber\\
&&+ P(\overline{Y},
\widehat{P}(\overline{Z},\overline{X}))\overline{W}
-R(T(\overline{X},\overline{Z}),\overline{Y})\overline{W}=0.
\end{eqnarray}
Setting $\overline{W}=\overline{\eta}$ into the above relation,
noting that $K \circ \gamma = id_{\cp}$,   $K\circ \beta=0$ and
$\widehat{S}=0$, it follows that
\begin{eqnarray*}
 &&(\nabla_{\gamma\overline{X}}\widehat{R})(\overline{Y},\overline{Z})
 -{R}(\overline{Y},\overline{Z})\overline{X}
   + (\nabla_{\beta\overline{Y}}\widehat{P})(\overline{Z},\overline{X})-
   (\nabla_{\beta
   \overline{Z}}\widehat{P})(\overline{Y},\overline{X})\nonumber\\
&&- \widehat{P}(\overline{Z},\widehat{P}(\overline{Y},\overline{X}))
+\widehat{R}(T(\overline{X},\overline{Y}),\overline{Z})+
\widehat{P}(\overline{Y}, \widehat{P}(\overline{Z},\overline{X}))
-\widehat{R}(T(\overline{X},\overline{Z}),\overline{Y})=0.
\end{eqnarray*}
Applying the  cyclic sum
$\mathfrak{S}_{\overline{X},\overline{Y},\overline{Z}}$ on the above
equation, taking  \textbf{(a)} into account, we get
\begin{equation}\label{eq.8}
    \mathfrak{S}_{\overline{X},\overline{Y},\overline{Z}}(
    \nabla_{\gamma\overline{X}}\widehat{R})(\overline{Y},\overline{Z})=0.
\end{equation}
Substituting (\ref{eq.1}) into (\ref{eq.8}), using
$(\nabla_{\gamma\overline{X}}\ell)(\overline{Y})=L^{-1}\hbar(\overline{X},\overline{Y})$
(Lemma \ref{ell}(b)), we have
\begin{eqnarray*}
   &&
    L(\nabla_{\gamma\overline{Z}}k)\set{(\ell(\overline{X})\overline{Y}-\ell(\overline{Y})\overline{X})} +
    L(\nabla_{\gamma\overline{Y}}k)\set{(\ell(\overline{Z})\overline{X}-\ell(\overline{X})\overline{Z})} \\
    &&+L(\nabla_{\gamma\overline{X}}k)\set{(\ell(\overline{Y})\overline{Z}-\ell(\overline{Z})\overline{Y})}
    +k\ell(\overline{Z})\set{(\ell(\overline{X})\overline{Y}-\ell(\overline{Y})\overline{X})}\\
   &&+k\ell(\overline{Y})\set{(\ell(\overline{Z})\overline{X}-\ell(\overline{X})\overline{Z})}
   +k\ell(\overline{X})\set{(\ell(\overline{Y})\overline{Z}-\ell(\overline{Z})\overline{Y})}\\
   &&+kL\set{(\hbar(\overline{X},\overline{Z})\overline{Y}-\hbar(\overline{Y},\overline{Z})\overline{X})}
   +kL\set{(\hbar(\overline{Z},\overline{Y})\overline{X}-\hbar(\overline{X},\overline{Y})\overline{Z})}\\
   &&+kL\set{(\hbar(\overline{Y},\overline{X})\overline{Z}-\hbar(\overline{Z},\overline{X})\overline{Y})}=0.
\end{eqnarray*}
Setting $\overline{Z}=\overline{\eta}$ into the above relation,
noting that $\ell(\overline{\eta})=L$ ,$\hbar(\overline{\eta},.)=0$
(Lemma \ref{ell}(c))
 and $\nabla_{\gamma\overline{\eta}}k=0$, we
 conclude that
 \begin{equation}\label{eq.2a}
   L^{2}\set{\nabla_{\gamma\overline{X}}k\,\phi(\overline{Y})
   -\nabla_{\gamma\overline{Y}}k\, \phi(\overline{X})}=0,
 \end{equation}
where $\phi$ is a vector $\pi$-form defined by
\begin{equation}\label{def.phi}
    g(\phi(\overline{X}),\overline{Y}):=\hbar(\overline{X},\overline{Y}).
\end{equation}

Taking the trace of both sides of (\ref{eq.2a}) with respect to
$\overline{Y}$, noting that $Tr(\phi)=n-1$ \cite{r86}, it follows
that
\begin{equation*}
(n-2)\nabla_{\gamma \overline{X}}k=0.
\end{equation*}
Consequently,
\begin{equation}\label{eq.9}
\nabla_{\gamma \overline{X}}k=0 \text{ \ \ for \ all  \ }
\overline{X}\in\cp, \text{ \ if\ }n\geq3.
\end{equation}
\par
 Now, applying the $v$-covariant derivative with respect to $\overline{Y}$
on both sides of (\ref{eq.7}), yields
\begin{equation*}
  \ell(\overline{Y})
  \nabla_{\beta\overline{X}}k+L(\stackrel{2}\nabla\stackrel{1}\nabla
  k)(\overline{X},\overline{Y})=
  L^{-1}\hbar(\overline{X},\overline{Y})(\nabla_{\beta\overline{\eta}}k)+\ell(\overline{X})
  (\stackrel{2}\nabla\stackrel{1}\nabla k)(\overline{\eta},\overline{Y}).
\end{equation*}
Since, $\stackrel{2}\nabla\stackrel{1}\nabla
  k=\stackrel{1}\nabla\stackrel{2}\nabla
  k=0$ (Lemma \ref{ricc} and (\ref{eq.9})), the above relation reduces to (provided that $n\geq3$)
\begin{equation*}
  \ell(\overline{Y})
  \nabla_{\beta\overline{X}}k=
  L^{-1}\hbar(\overline{X},\overline{Y})(\nabla_{\beta\overline{\eta}}k).
\end{equation*}
Setting $\overline{Y}=\overline{\eta}$ into the above equation,
noting that $\ell(\overline{\eta})=L$ and
$\hbar(.,\overline{\eta})=0$, it follows that
$\nabla_{\beta\overline{\eta}}k=0$. Consequently, again by
(\ref{eq.7}),
\begin{equation}\label{eq.11}
\nabla_{\beta \overline{X}}k=0\text{ \ \ for \ all  \ }
\overline{X}\in\cp, \text{ \ if\ }n\geq3.
\end{equation}
Now, Equations  (\ref{eq.9}) and (\ref{eq.11}) imply  that $k$ is a
constant if $n\geq3$.\\
 This complete the proof.
\end{proof}

\bigskip

\begin{thm}\label{th.2}Let $(M,L)$ be a Finsler manifold with
dimension $n\geq3$ and let $q\neq 1$ be an arbitrary real number. If
the $h$-curvature tensor $R$ satisfies
\begin{equation}\label{eq.12}
    R(\overline{X},\overline{Y})\overline{Z}=k\,\mathfrak{A}_{\overline{X},\overline{Y}}
    \set{q\,g(\overline{X},\overline{Z})\overline{Y}+(q-2)\set{L^{-1}g(\overline{X},\overline{Z})
    \ell(\overline{Y})\overline{\eta}-\ell(\overline{Y})\ell(\overline{Z})\overline{X}}},
    \end{equation}
where $k(x,y)$ is  a homogenous function of degree $0$ in $y$, then
\begin{description}
   \item[(a)] $k$ is a constant.
    \item[(b)] If $k\neq0$, we have:
    \begin{description}
     \item[1)]
      $P(\overline{X},\overline{Y})\overline{Z}=P(\overline{Y},\overline{X})\overline{Z}$
      \emph{(}i.e.,\,$(M,L)$ is symmetric\emph{)},
     \item[2)]
      $S(\overline{X},\overline{Y})\overline{Z}=\displaystyle{\frac{2-q}{2(q-1)L^2}}\set
      {\hbar(\overline{X},\overline{Z})\phi(\overline{Y})-\hbar(\overline{Y},\overline{Z})\phi(\overline{X})}$.
    \end{description}
\end{description}
\end{thm}
\begin{proof} ~\par
 \noindent\textbf{(a)} Setting $\overline{Z}=\overline{\eta}$ into
 (\ref{eq.12}), we get
\begin{equation}\label{eq.13}
    \widehat{R}(\overline{X},\overline{Y})=2(q-1)kL
      \set{\ell(\overline{X})\overline{Y}-\ell(\overline{Y})\overline{X}}.
\end{equation}
From which, together with Theorem \ref{prop.1}, the result follows.

\vspace{5pt}
 \noindent\textbf{(b)\,1).} Applying  the $v$-covariant derivative with
 respect to $\overline{W}$ on both sides of (\ref{eq.12}),
 we get
\begin{equation*}
    (\nabla_{\beta
    \overline{W}}R)(\overline{X},\overline{Y},\overline{Z})=0.
\end{equation*}
From which, together (\ref{eq.4}), it follows that
\begin{equation}\label{eq.14}
   \mathfrak{S}_{\overline{X},\overline{Y},\overline{Z}}\,
P(\widehat{R}(\overline{X},\overline{Y}),
\overline{Z})\overline{W}=0.
\end{equation}
In view of   (\ref{eq.13}), noting that $k\neq0$, (\ref{eq.14})
implies that
\begin{eqnarray*}
  &&2(q-1)L\set{P(\ell(\overline{X})\overline{Y}-\ell(\overline{Y})\overline{X},
\overline{Z})\overline{W}}+2(q-1)L\set{P(\ell(\overline{Y})\overline{Z}-\ell(\overline{Z})\overline{Y},
\overline{X})\overline{W}}\\
  &&+2(q-1)L\set{P(\ell(\overline{Z})\overline{X}-\ell(\overline{X})\overline{Z},
\overline{Y})\overline{W}}=0.
\end{eqnarray*}
 Setting $\overline{Z}=\overline{\eta}$ into the above equation, taking into account
 the fact that $\ell(\overline{\eta})=L$ and  $P(.,\overline{\eta}).=P(\overline{\eta},.).=0$ \cite{r96}, we get
 \begin{equation*}
2(q-1)L\set{P(\overline{X},
\overline{Y})\overline{W}-P(\overline{Y},
\overline{X})\overline{W}}=0.
 \end{equation*}
Hence, the result follows.

\vspace{5pt}
 \noindent\textbf{(b)\,2).} Taking the cyclic sum
 $\mathfrak{S}_{\overline{X},\overline{Y},\overline{Z}}$ of
  (\ref{eq.r1a}) and using \textbf{(b)}1), we obtain
\begin{eqnarray}\label{eq.r1b}
 \mathfrak{S}_{\overline{X},\overline{Y},\overline{Z}}\set{
 (\nabla_{\gamma\overline{X}}R)(\overline{Y},\overline{Z},\overline{W})-
S(\widehat{R}(\overline{Y},\overline{Z}),\overline{X})\overline{W}}=0.
\end{eqnarray}
\par
 On the other hand, by taking the $v$-covariant derivative of both
 sides of (\ref{eq.12}), using $(\nabla_{\gamma\overline{X}}L)=\ell(\overline{X})$,
  $(\nabla_{\gamma\overline{X}}\ell)(\overline{Y})=L^{-1}\hbar(\overline{X},\overline{Y})$ and
 $\nabla_{\gamma\overline{X}}g=0$, we get
\begin{eqnarray*}
  (\nabla_{\gamma\overline{X}}R)(\overline{Y},\overline{Z},\overline{W}) &=&
  k(q-2)\mathfrak{A}_{\overline{X},\overline{Y}}\{g(\overline{X},\overline{W})
  \hbar(\overline{Z},\overline{Y})\frac{\overline{\eta}}{L^{2}} +
  g(\overline{X},\overline{W})
  \ell(\overline{Y})\frac{\phi(\overline{Z})}{L}\\
 && -\hbar({\overline{Z},\overline{Y}})\ell(\overline{W})\frac{\overline{X}}{L}
  -\hbar(\overline{Z},\overline{W})\ell(\overline{Y})\frac{\overline{X}}{L}\}.
\end{eqnarray*}
Taking the cyclic sum
 $\mathfrak{S}_{\overline{X},\overline{Y},\overline{Z}}$ of both
 sides of the above equation and then setting
 $\overline{Z}=\overline{\eta}$, it follows that
\begin{equation}\label{eq.16}
\mathfrak{S}_{\overline{X},\overline{Y},\,\overline{\eta}}
(\nabla_{\gamma\overline{X}}R)(\overline{Y},\overline{\eta},\overline{W})=2k(q-2)
\set{\hbar(\overline{Y},\overline{W})\phi(\overline{X})-
\hbar(\overline{X},\overline{W})\phi(\overline{Y})}.
\end{equation}
 In view of (\ref{eq.13}), noting that $S(.,\overline{\eta}).=0$ and
 $S$ is antisymmetric \cite{r96}, we obtain
\begin{equation}\label{eq.17}
\mathfrak{S}_{\overline{X},\overline{Y},\,\overline{\eta}}
S(\widehat{R}(\overline{Y},\overline{\eta}),\overline{X})\overline{W}
=4kL^{2}(q-1)S(\overline{X},\overline{Y})\overline{W}.
\end{equation}
Therefore, by setting $\overline{Z}=\overline{\eta}$ into
(\ref{eq.r1b}), taking (\ref{eq.16}) and (\ref{eq.17}) into account,
the result follows.
\end{proof}

\begin{cor} Akbar-Zadeh's theorem \emph{\cite{r58a}} is a special case of Theorem
\ref{th.2}, for which $q=2$.
\end{cor}

\begin{cor} If
the $h$-curvature tensor $R$ of $(M,L)$, ($n\geq3$), satisfies
\begin{equation*}
    R(\overline{X},\overline{Y})\overline{Z}=k\mathfrak{A}_{\overline{X},\overline{Y}}
    \set{\{g(\overline{X},\overline{Z})\frac{\overline{\eta}}{L}
    -\ell(\overline{Z})\overline{X}\}\ell(\overline{Y})},
\end{equation*}
 then, $k$ is a constant and,  moreover, if $k\neq0$, we have:
\begin{description}
   \item[(a)]$(M,L)$ is symmetric.
    \item[(b)] $S(\overline{X},\overline{Y})\overline{Z}=\displaystyle{\frac{-1}{L^{2}}}\set
      {\hbar(\overline{X},\overline{Z})\phi(\overline{Y})-\hbar(\overline{Y},\overline{Z})\phi(\overline{X})}.$
 \end{description}
\end{cor}

\Section{Second generalization of Akbar-Zadeh's theorem}

 In this section, we give a second intrinsic  generalization of
 Akbar-Zadeh's
 theorem.

 \begin{thm}\label{thm.3} If the $h$-curvature tensor $R$ of $(M,L)$, $\emph{dim} M \geq3$,
 satisfies
\begin{equation}\label{eq.18}
    R(\overline{X},\overline{Y})\overline{Z}=k
    \set{g(\overline{X},\overline{Z})\overline{Y}-g(\overline{Y},\overline{Z})\overline{X}
    +\omega(\overline{X},\overline{Y})\overline{Z}},
\end{equation}
 where  $\omega$ is an indicatory antisymmetric h(0) $\pi$-tensor field of type
 $(1,3)$ and $k(x,y)$ is an h(0)-function on $TM$, then
\begin{description}
   \item[(a)] $k$ is a constant.
    \item[(b)] If $k\neq0$, we have:
\begin{description}
  \item[1)]
  $P(\overline{X},\overline{Y})\overline{Z}-P(\overline{X},\overline{Y})\overline{Z}=
  L^{-2}(\nabla_{\beta\overline{\eta}}\omega)(\overline{X},\overline{Y},\overline{Z}).$
 \item[2)] $S(\overline{X},\overline{Y})\overline{Z}=\displaystyle{\frac{1}{L^{2}}}
  \set{\displaystyle{\frac{1}{2kL^{2}}}(\stackrel{1}\nabla\stackrel{1}\nabla\omega)
  (\overline{\eta},\overline{\eta}, \overline{X}, \overline{Y},
  \overline{Z})+\omega(\overline{X},\overline{Y})\overline{Z}}.$
\end{description}
\end{description}
 \end{thm}

 \begin{proof}~\par
 \noindent\textbf{(a)} Follows from Therorem  \ref{prop.1} by setting
 $\overline{Z}=\overline{\eta}$ into (\ref{eq.18}).

\vspace{5pt}
 \noindent\textbf{(b)\,1)}. By (\ref{eq.18}), we have
 \begin{equation}\label{eq.19}
    \widehat{R}(\overline{X},\overline{Y})=kL
    \set{\ell(\overline{X})\overline{Y}-\ell(\overline{Y})\overline{X}},
\end{equation}
and by (\ref{eq.4}), we have
\begin{equation}\label{eq.24}
   \mathfrak{S}_{\overline{X},\overline{Y},\overline{\eta}}\,
\{(\nabla_{\beta \overline{X}}R)(\overline{Y},
\overline{\eta},\overline{W})+P(\widehat{R}(\overline{X},\overline{Y}),
\overline{\eta})\overline{W}\}=0.
\end{equation}
 Now, substituting (\ref{eq.18}) and (\ref{eq.19}) into
(\ref{eq.24}), we obtain
\begin{equation*}
k\set{(\nabla_{\beta\overline{\eta}}\omega)(\overline{X},\overline{Y},\overline{W})-
L^{2}\set{P(\overline{X},\overline{Y})\overline{Z}-P(\overline{Y},\overline{X})\overline{Z}}}=0
\end{equation*}
From which, since $k\neq0$, the result follows.

\vspace{5pt}
 \noindent\textbf{(b)\,2)}. Taking the cyclic sum
 $\mathfrak{S}_{\overline{X},\overline{Y},\overline{Z}}$ of  (\ref{eq.r1a}), we obtain
\begin{eqnarray}\label{eq.r1c}
 &&\mathfrak{S}_{\overline{X},\overline{Y},\overline{Z}}\{
 (\nabla_{\gamma\overline{X}}R)(\overline{Y},\overline{Z},\overline{W})+
 (\nabla_{\beta\overline{Y}}P)(\overline{Z},\overline{X},\overline{W})\nonumber\\
 &&-
   (\nabla_{\beta
   \overline{Z}}P)(\overline{Y},\overline{X},\overline{W})-
S(\widehat{R}(\overline{Y},\overline{Z}),\overline{X})\overline{W}\}=0.
\end{eqnarray}
In view of \textbf{1)} above, it follows that
\begin{equation*}
    (\nabla_{\beta\overline{W}}P)(\overline{X},\overline{Y},\overline{Z})
    - (\nabla_{\beta\overline{W}}P)(\overline{Y},\overline{X},\overline{Z})=
  L^{-2}(\stackrel{1}\nabla\stackrel{1}\nabla\omega)
  (\overline{W},\overline{\eta}, \overline{X}, \overline{Y},
  \overline{Z})
\end{equation*}
 From which, we get
\begin{equation}\label{eq.a1}
    \mathfrak{S}_{\overline{X},\overline{Y},\overline{\eta}}\set{
    (\nabla_{\beta\overline{\eta}}P)(\overline{X},\overline{Y},\overline{Z})
    - (\nabla_{\beta\overline{\eta}}P)(\overline{Y},\overline{X},\overline{Z})}=
  L^{-2}(\stackrel{1}\nabla\stackrel{1}\nabla\omega)
  (\overline{\eta},\overline{\eta}, \overline{X}, \overline{Y},
  \overline{Z}).
\end{equation}
\par
On the other hand, noting that $\omega$ is homogenous of degree
zero, we obtain
\begin{eqnarray}
   \mathfrak{S}_{\overline{X},\overline{Y},\overline{\eta}}
 (\nabla_{\gamma\overline{X}}R)(\overline{Y},\overline{\eta},\overline{W})&=&
 (\nabla_{\gamma\overline{X}}R)(\overline{Y},\overline{\eta},\overline{W})
 +(\nabla_{\gamma\overline{Y}}R)(\overline{\eta},\overline{X},\overline{W})\nonumber \\
 &&+(\nabla_{\gamma\overline{\eta}}R)(\overline{X},\overline{Y},\overline{W})\nonumber \\
   &=&  2k\omega(\overline{X},\overline{Y})\overline{W}.\label{eq.111}\\
 \mathfrak{S}_{\overline{X},\overline{Y},\overline{\eta}}
S(\widehat{R}(\overline{X},\overline{Y}),\overline{\eta})\overline{W}&=&
   2kL^{2}S(\overline{X},\overline{Y})\overline{W}.\label{eq.112}
   \end{eqnarray}
Setting $\overline{Z}=\overline{\eta}$ into (\ref{eq.r1c}), taking
into account (\ref{eq.a1}), (\ref{eq.111})  and (\ref{eq.112}), the
result follows.
\end{proof}

\begin{cor} Akbar-Zadeh's theorem {\emph{\cite{r58a}}} is obtained from the above Theorem
by  letting $\omega=0$.
\end{cor}

\begin{cor} A Finsler manifold $(M,L)$ is $S_{3}$-like if $\omega$ in
Theorem \ref{thm.3} is given by
\begin{equation}\label{eq.a2}
    \omega(\overline{X},\overline{Y})\overline{Z}=S\set{\hbar(\overline{X},\overline{Z}
    )\phi(\overline{Y})-\hbar(\overline{Y},\overline{Z})\phi(\overline{X})},
\end{equation}
where $\phi$ is given by (\ref{def.phi}) and  $S(x)$ is a scalar
function independent of $y$.
\end{cor}
\begin{proof}
From Theorem \ref{thm.3}(b) and (\ref{eq.a2}), the $v$-curvature
tensor $S$ takes the form:
\begin{equation*}
    S(\overline{X},\overline{Y})\overline{Z}=\frac{1}{L^{2}}
  \set{S+\frac{(\stackrel{1}\nabla\stackrel{1}\nabla S)
  (\overline{\eta},\overline{\eta})}{2kL^{2}}}
  \set{\hbar(\overline{X},\overline{Z}
    )\phi(\overline{Y})-\hbar(\overline{Y},\overline{Z})\phi(\overline{X})}.
\end{equation*}
As  the $v$-curvature tensor $S$ is written in the above  form, then
the term  $$\set{S+\frac{(\stackrel{1}\nabla\stackrel{1}\nabla S)
  (\overline{\eta},\overline{\eta})}{2kL^{2}}}$$ depends on $x$ only \cite{r44a},
  and so $(M,L)$ is $S_{3}$-like.
\end{proof}

\begin{cor} If the scalar function $S(x)$ in (\ref{eq.a2}) is  constant, we
have:
\begin{description}
  \item[(a)]$P(\overline{X},\overline{Y})\overline{Z}=P(\overline{Y},\overline{X})\overline{Z}$.
  \item[(b)] $S(\overline{X},\overline{Y})\overline{Z}=\frac{S}{L^{2}}\set{\hbar(\overline{X},\overline{Z}
    )\phi(\overline{Y})-\hbar(\overline{Y},\overline{Z})\phi(\overline{X})}$.
\end{description}
\end{cor}

\begin{cor} If the tensor field $\omega$ in Theorem \ref{thm.3} is
given by
 $$\omega(\overline{X},\overline{Y})\overline{Z}=
 \mathfrak{A}_{\overline{X},\overline{Y}}\set{H(\overline{X},\overline{Z}
    )\phi(\overline{Y})+\hbar(\overline{X},\overline{Z})H_{o}(\overline{Y})},$$
where $H$ is a symmetric indicatory  h(0) 2-scalar $\pi$-form and
$H(\overline{X},\overline{Y})=:g(H_{o}(\overline{X}),\overline{Y})$,
then $(M,L)$ is $S_{4}$-like, that is,
 $$S(\overline{X},\overline{Y})\overline{Z}=\frac{1}{L^{2}}
 \mathfrak{A}_{\overline{X},\overline{Y}}\set{\mu(\overline{X},\overline{Z}
    )\phi(\overline{Y})+\hbar(\overline{X},\overline{Z})\mu_{o}(\overline{Y})},$$
where\,
$\mu(\overline{X},\overline{Y})=\set{H(\overline{X},\overline{Y})+\frac{(
\stackrel{1}\nabla\stackrel{1}\nabla H)
  (\overline{\eta},\overline{\eta},\overline{X},\overline{Y})}{2kL^{2}}}$.
\end{cor}

\begin{proof} The proof is clear and we omit it.
\end{proof}

\bigskip
\noindent\textbf{Concluding remark.}\emph{ It should be noted that
the outcome of this work is twofold. Firstly,   the local
expressions of the obtained results, when calculated, coincide with
the existing local results (\cite{r58a} , \cite{hojo}). Secondly,
new global proofs have been established.}

\bigskip
\bigskip
\noindent\textbf{Acknowledgments}

\vspace{5pt}
 The author express his sincere
thanks to referees and my Professor Nabil L. Youssef for their
valuable suggestions and comments.

\providecommand{\bysame}{\leavevmode\hbox
to3em{\hrulefill}\thinspace}
\providecommand{\MR}{\relax\ifhmode\unskip\space\fi MR }
\providecommand{\MRhref}[2]{%
  \href{http://www.ams.org/mathscinet-getitem?mr=#1}{#2}
} \providecommand{\href}[2]{#2}

\end{document}